\font\Sc=cmcsc10
\font\Ch=msbm9

\font\small=cmr8

\font\ninerm=cmr9

\parindent=0pt
\def\nl{\bigskip\noindent}
\newcount\subsectionno
\def\resetss{\subsectionno=0}
\def\subsection#1
{\global\advance\subsectionno by 1
 \goodbreak\par\bigskip\bigskip\noindent{\Sc \the\sectionno.\the\subsectionno \hskip .1in #1}
\par\bigskip
}
\newcount\sectionno
\sectionno = 0
\def\section#1
{\global\advance\sectionno by 1 \resetss
 \goodbreak\par\bigskip\bigskip\noindent{\Sc \the\sectionno. #1}
\par
}
\newcount\refctr  \refctr = 0
\def\enumerate#1 #2\endenumerate{\advance\refctr by 1
  \expandafter\edef\csname #1\endcsname{\the\refctr}%
  \if !#2! \else\enumerate #2\endenumerate\fi}
\newcount\numerobib \numerobib=0
\def\bibitem{\advance\numerobib by 1
    \par\noindent\item{[\number\numerobib]}}

\def\Q{\hbox{\Ch Q}}

\newcount\elemno
\elemno=0
\def\numero{\global\advance\elemno by 1
            {\the\elemno}.}
\def\TITRECO#1{\goodbreak\nl {\bf #1}\ {\Sc \numero}\ \ \ \ }

\def\Ex{\TITRECO{Example}}
\outer\def\Thm#1\par{\bigbreak
         \noindent{{\bf Theorem}\ {\Sc \numero\quad}}{\sl#1}\par
          \ifdim\lastskip<\medskipamount \removelastskip\penalty55\medskip\fi}
\outer\def\Prop#1\par{\bigbreak
         \noindent{{\bf Proposition}\ {\Sc \numero\quad}}{\sl#1}\par
          \ifdim\lastskip<\medskipamount \removelastskip\penalty55\medskip\fi}
\outer\def\Cor#1\par{\bigbreak
         \noindent{{\bf Corollary}\ {\Sc \numero\quad}}{\sl#1}\par
          \ifdim\lastskip<\medskipamount \removelastskip\penalty55\medskip\fi}
\outer\def\Lem#1\par{\bigbreak
         \noindent{{\bf Lemma}\ {\Sc \numero\quad}}{\sl#1}\par
          \ifdim\lastskip<\medskipamount \removelastskip\penalty55\medskip\fi}
\outer\def\Conj#1\par{\bigbreak
         \noindent{{\bf Conjecture}\ {\Sc \numero\quad}}{\sl#1}\par
          \ifdim\lastskip<\medskipamount \removelastskip\penalty55\medskip\fi}
\def\Proof{\nl {\sl Proof\/}.\ \ \ }
\newcount\noeqn
\noeqn = 0
\def\mark{\global\advance\noeqn by 1
\eqno (\the\noeqn)}
\newdimen\Squaresize \Squaresize=14pt
\newdimen\Thickness \Thickness=0.5pt

\def\Square#1{\hbox{\vrule width \Thickness
   \vbox to \Squaresize{\hrule height \Thickness\vss
      \hbox to \Squaresize{\hss#1\hss}
   \vss\hrule height\Thickness}
\unskip\vrule width \Thickness}
\kern-\Thickness}

\def\Vsquare#1{\vbox{\Square{$#1$}}\kern-\Thickness}
\def\Blk{\omit\hskip\Squaresize}

\def\Young#1{
\vbox{\smallskip\offinterlineskip
\halign{&\Vsquare{##}\cr #1}}}

\newdimen\squaresize \squaresize=3pt
\newdimen\thickness \thickness=0.2pt

\def\square#1{\hbox{\vrule width \thickness
   \vbox to \squaresize{\hrule height \thickness\vss
      \hbox to \squaresize{\hss#1\hss}
   \vss\hrule height\thickness}
\unskip\vrule width \thickness}
\kern-\thickness}

\def\vsquare#1{\vbox{\square{$#1$}}\kern-\thickness}

\def\thisbox#1{\kern-.09ex\fbox{#1}}
\def\downbox#1{\lower1.200em\hbox{#1}}
\def\picture #1 by #2 (#3){
  \vbox to #2{
    \hrule width #1 height 0pt depth 0pt
    \vfill
    \special{picture #3} 
    }
  }

\def\scaledpicture #1 by #2 (#3 scaled #4){{
  \dimen0=#1 \dimen1=#2
  \divide\dimen0 by 1000 \multiply\dimen0 by #4
  \divide\dimen1 by 1000 \multiply\dimen1 by #4
  \picture \dimen0 by \dimen1 (#3 scaled #4)}
  }

\overfullrule=0pt
\baselineskip 14pt
\settabs 8 \columns
\+ \hfill&\hfill \hfill & \hfill \hfill   & \hfill \hfill  & \hfill \hfill   &
 & \hfill \hfill\cr
\parindent=.5truein
\hfuzz=3.44182pt
\hsize 6truein
\font\small=cmr6
\font\title=cmbx10 scaled\magstep2
\font\normal=cmr12 
\font\small=cmr6
\font\ninerm=cmr9 

\font\bol=cmbx12
\normal

\vsize=8truein
\def\today{\ifcase\month\or
January\or February\or March\or April\or may\or June\or
July\or August\or September\or October\or November\or
December\fi
\space\number\day, \number\year}
\def\today{\ifcase\month\or
January\or February\or March\or April\or May\or June\or
July\or August\or September\or October\or November\or
December\fi
\space\number\day, \number\year}
\headline={\small \hfill Ribbons and Hall-Littlewood symmetric functions
\hfill \today \hskip .5in Page \folio} \footline={\hfil}

\centerline{\bol Ribbon Operators and}
\centerline{\bol Hall-Littlewood Symmetric Functions}
\font\Ch=msbm9

\vskip .2in
\centerline{\bf Mike Zabrocki}
\centerline{\it{ LaCIM}}
\centerline{\it{ Universit\'e du Qu\'ebec \`a Montr\'eal}}
\centerline{\it{ Montr\'eal (Qu\'ebec) H3C 3P8}}
\centerline{ \tt zabrocki@math.uqam.ca}

\def\la{{\lambda}}
\def\coeff{{\Big|}}
\def\ht #1{{\overline{#1}}}

\def\endofproof {\hskip .1in $\diamondsuit$}
\def\keq{\simeq_k}

\def\Sw {{\tilde S}}

\def\X {\hbox{\bf X}}

\def\pointir{\discretionary{.}{}{.\kern.35em---\kern.7em}\nobreak
\hskip 0em plus .3em minus .4em }
\def\article#1|#2|#3|#4|#5|#6|#7|
    {{\leftskip=7mm\noindent
     \hangindent=0mm\hangafter=1
     \llap{#1.\hskip.35em}{#2,}
     #3, {\it #4} \nobreak {\bf #5} \nobreak {#6},
     \nobreak\ #7.\par}}
\def\unarticle#1|#2|#3|#4|
    {{\leftskip=7mm\noindent
     \hangindent=0mm\hangafter=1
     \llap{#1.\hskip.35em}{#2,}
     #3, {\it #4}, \nobreak to appear.\par}}
\def\livre#1|#2|#3|#4|
    {{\leftskip=7mm\noindent
    \hangindent=0mm\hangafter=1
    \llap{#1.\hskip.35em}{#2,}
    {``#3,''} #4.\par}}
\vskip .3in

\noindent
{\Sc Abstract.} {\ninerm Given a partition $\la = (\la_1, \la_2, \ldots
\la_k)$, let $\la^{rc} = (\la_2-1, \la_3-1, \ldots
\la_k-1)$.  It is easily seen that the diagram $\la\slash \la^{rc}$ is connected
and has no $2 \times 2$ subdiagrams which we shall refer to as a ribbon.  To each
ribbon $R$, we associate a symmetric function operator $S^R$. 
We may define the major index of a ribbon $maj(R)$
to be the major index of any permutation that fits the ribbon.  This paper
is concerned with the operator $H_{1^k}^q = \sum_R q^{maj(R)} S^R$ where the
sum is over all $2^{k-1}$ ribbons of size $k$.  We show here that
$H_{1^k}^q$ has truly remarkable properties, in particular
that it is a Rodriguez operator that adds a column to the Hall-Littlewood
symmetric functions.  We believe that some of the tools we introduce
here to prove our results should also be of independent interest and may
be useful to establish further symmetric function identities. }

\section{Introduction}

The Schur functions indexed by a sequence of integers $(p_1, p_2, \ldots, p_k)$
can be defined by the Jacobi-Trudi identity 
$s_{(p_1, p_2, \ldots, p_k)} = det|h_{p_i + i-j} |_{1 \leq i,j \leq k}$.  
It is well known and easy to show that we have the
relation $$s_{(p_1, \ldots,p_i, p_{i+1},\ldots, p_k)}
= -s_{(p_1, \ldots,p_{i+1}-1, p_{i}+1,\ldots, p_k)}.\mark$$  

Let us recall that for a given symmetric function $f$ it is customary to
denote by $f^\perp$ the operator that is dual to multiplication by $f$
with respect to the Hall inner product.  We shall make crucial use here of
the Bernstein [5] operator

$$S_m= \sum_{k\geq 0} (-1)^k h_{m+k} e_k^\perp \mark$$
Its action on the Schur basis may be easily computed with
the formula $S_m s_{(p_1, p_2, \ldots, p_k)} = 
s_{(m, p_1, p_2, \ldots, p_k)}$ and relation $(1)$.

In particular we have the ``Rodriguez'' formula for Schur functions
indexed by a partition $\la$,

$$S_{\la_1} S_{\la_2} \cdots S_{\la_k} 1 = s_\la\mark$$

For a partition $\la = (\la_1 \geq \la_2 \geq \ldots \geq \la_m > 0)$
we set $\ell(\la) = m$.
Partitions here are drawn by the
French convention with the smallest part on the top.

A ribbon is a connected skew partition that contains no 
$2 \times 2$ blocks.
If $\la$ is a partition of length $k$, we set
$\lambda^{rc} = (\la_2-1,\la_3 - 1, \ldots, \la_k-1)$ 
(the partition with the first
row and column removed).  It is easy to see that
every ribbon partition will be 
$\la \slash \la^{rc}$ for
some partition $\la$.

Label the cells
in a ribbon diagram with the numbers $\{1, 2, \ldots, k-1 \}$
from left to right, top to bottom. 
This done, we let
$$D(R) = \{ i \in [1,k-1] : i+1^{st} \hbox{\ cell of $R$ lies below
the $i^{th}$ cell } \}\mark$$
and refer to it as the descent set of $R$.  The ribbons are therefore
in one to one correspondence with the subsets of $\{ 1, 2, \dots k-1\}$.
If $R$ is a ribbon of size $k$, we use $\ht{R}$ to denote the ribbon
whose descent set is $\{ 1, 2, \ldots k-1 \} - D(R)$.

We also let the `major'
index of $R$ be defined as

$$maj(R) = \sum_{i \in D(R)} i \mark$$
and we set $comaj(R) = {k \choose 2} - maj(R)$ to be the complementary
statistic.  Clearly we have the relation $maj(\ht{R}) = comaj(R)$.

Below we have listed all of the ribbon partitions of size $4$ with the
corresponding descent set which is a subset of $\{1,2,3\}$.  

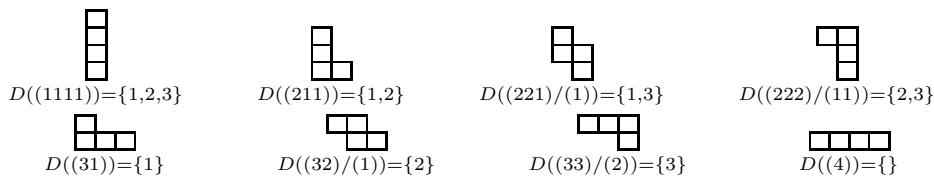
\Squaresize=7pt
\Thickness=0.5pt
\centerline{
$\Young{\cr\cr\cr\cr} \atop D({(1111)}) = \{1,2,3\}$ \hskip .3in
$\Young{\Blk\cr\cr\cr  & \cr} \atop D({(211)})=\{1,2\}$ \hskip .3in
$\Young{\Blk\cr\cr&\cr \Blk & \cr} \atop D({(221)\slash(1)}) = \{1,3\}$ \hskip .3in
$\Young{\Blk\cr&\cr\Blk&\cr \Blk & \cr} \atop D({(222)\slash(11)})=\{2,3\}$}

\centerline{
$\Young{\Blk\cr\Blk\cr\cr  & &\cr} \atop D({(31)})=\{1\}$ \hskip .5in
$\Young{\Blk\cr\Blk\cr&\cr \Blk & &\cr} \atop D({(32)\slash(1)})=\{2\}$ \hskip .4in
$\Young{\Blk\cr\Blk\cr & &\cr \Blk &\Blk & \cr} \atop D({(33)\slash(2)})=\{3\}$ \hskip .5in
$\Young{\Blk\cr\Blk\cr\Blk\cr & & & \cr} \atop D({(4)})=\{\}$}

\vskip .1in
We will use the symbol $R$ to represent an arbitrary ribbon and
the notation $R \models k$ to indicate that $R$ is a ribbon of size $k$.

For each ribbon $R$ of size $k$ we create an operator that raises the degree
of the symmetric function by $k$.  If $R = \la \slash \la^{rc}$ then set
$$S^R = s_{\la^{rc}}^\perp \Sw_{\la_1'} \Sw_{\la_2'} \cdots \Sw_{\la_{\la_{1}}'}.
\mark$$
where $\la_i'$ is the length of the $i^{th}$ column in the partition $\la$.

The action of this operator is very combinatorial in nature,
we attach a ribbon on the left of the Schur function and
reduce using the commutation relations $\Sw_a \Sw_b = -\Sw_{b-1} \Sw_{a+1}$
and $\Sw_a \Sw_{a+1} = 0$,
followed by the Littlewood-Richardson rule.  We present one large example
below.

$$S^{(432)/(21)} (s_{(33221)}) = -s_{(665)} - s_{(764)}- s_{(755)}$$

$${S^{\Young{& \cr\Blk&&\cr\Blk&\Blk& & \cr}}} \Bigg( {\Young{\cr&\cr & \cr& & \cr
& & \cr}} \Bigg) 
= \Young{\Blk& \Blk&\Blk& \Blk&\cr
\Blk& \Blk&\Blk& \Blk&&\cr
&&\Blk& \Blk&&\cr
 \Blk &  & & \Blk&&&\cr
\Blk& \Blk& &  &&&\cr} 
= - \Young{&&&&&\cr\Blk&&&&&&\cr\Blk&\Blk&&&&&\cr}
= - \Young{&&&&\cr&&&&&\cr&&&&&\cr}
- \Young{&&&\cr&&&&&\cr&&&&&&\cr}
- \Young{&&&&\cr&&&&\cr&&&&&&\cr}$$

Define the symmetric functions $H_\la[X;q] = \sum_\mu K_{\la\mu}(q) s_\mu[X]$
where $K_{\la\mu}(q)$ is the Kostka-Folkes polynomial.
There is an operator $H_m^q$ that adds a row to this symmetric function
when $m$ is larger than $\la_1$
that is due to Jing [2]. In particular, it yields the `Rodriguez' formula

$$H_{\la_1}^q H_{\la_2}^q \cdots H_{\la_k}^q  1= 
H_{\la}[X;q]. \mark $$

Our main result here is the construction of an operator $H_{1^k}^q$
which adds a {\it column} to the partition indexing $H_{\la}[X;q]$.
More precisely, we show here that the operator

$$H_{1^k}^q = \sum_{R\models k }q^{comaj(R)} S^R \mark$$
has the following remarkable properties.

\Thm For all $k \geq 0$, $$H_{m+1}^q H_{1^k}^q 
= H_{1^{k+1}}^q H_{m}^q.\mark$$
As a result, we have for $\ell = \ell(\la)$
$$H_{1^{\la_1}}^q H_{1^{\la_2}}^q \cdots H_{1^{\la_\ell}}^q 1 
= H_{\la'}[X;q] \mark$$
where $\la'$ is the conjugate partition to $\la$.
\vskip .1in

The property that $H_{1^k}^q$ adds a column to the Hall-Littlewood
symmetric functions is a consequence of the commutation relation, since
$$H_{1^k}^q H_\la[X,q] = 
H_{\la_1+1}^q H_{\la_2+1}^q \cdots H_{\la_k+1}^q H_{1^0}^q(1) =
H_{(\la_1+1,\la_2+1,\cdots,\la_k+1)}[X;q].\mark $$
 
This result is the end product of a number of very interesting identities
satisfied by ribbon operators.  Our basic tool in establishing them is
a truly remarkable new involution in the theory of symmetric functions.
However we need to postpone the statement of these further results to the
next section after after we introduce some less familiar notation.

We should mention that $(10)$ is a rather surprising extension to the
general case of the classical identity
$$H_{1^k}[X;q] = \sum_{\sigma \in S_k} q^{comaj(\sigma)} 
s_{\la(\sigma)}[X] \mark$$
where the sum is over the symmetric group $S_k$ and $\la(\sigma)$ denotes
the shape of the standard tableaux corresponding to $\sigma$ under
Robinson-Schenstead correspondence.  In fact, by grouping terms according
to descent sets we derive from $(12)$ that
$$H_{1^k}[X;q] = \sum_{R \models k} q^{comaj(R)} 
s_{R}[X]. \mark$$

\Ex 
We see that for $k=3$, formula $(13)$ reduces

$$H_{(1^3)}[X;q] = \Young{\cr\cr\cr} + q \Young{&\cr\Blk&\cr}
+ q^2 \Young{\cr&\cr} + q^3 \Young{&&\cr}\mark$$

To compute the Hall-Littlewood symmetric function $H_{(2^3)}[X;q]$ we
act on the symmetric function $H_{(1^3)}[X;q] =
s_{(1^3)} + (q+q^2) s_{(21)} + q^3 s_{(3)}$ with each ribbon of size $3$.

$$\Young{\cr\cr\cr}H_{(1^3)}[X;q] = \Young{&\cr&\cr&\cr} 
+ (q+q^2) \Young{\cr&\cr&&\cr} + q^3 \Young{\cr\cr&&&\cr}\mark$$

$$\eqalign{q\Young{&\cr\Blk&\cr}H_{(1^3)}[X;q] &= q \Young{\Blk&\Blk&\cr&&\cr\Blk&&\cr} 
+ (q^2+q^3) \Young{&&\cr\Blk&&&\cr} + q^4 \Young{&\cr\Blk&&&&\cr}\cr
&= (q^2+q^3) \left(\Young{&\cr&&&\cr} + \Young{&&\cr&&\cr}\right)
+q^4 \left(\Young{&\cr&&&\cr} + \Young{\cr&&&&\cr}\right)}\mark$$

$$\eqalign{q^2\Young{&\Blk\cr&\cr}H_{(1^3)}[X;q] &= 
q^2 \Young{\Blk&\Blk&\cr&\Blk&\cr&&\cr} 
+ (q^3+q^4) \Young{&\Blk&\cr&&&\cr} + q^5 \Young{\cr&&&&\cr}\cr
&= -q^2 \Young{&&\cr&&\cr}
+q^5 \Young{\cr&&&&\cr}}\mark$$

$$\eqalign{q^3\Young{&&\cr}H_{(1^3)}[X;q] &= 
q^3 \Young{\Blk&\Blk&\Blk&\cr\Blk&\Blk&\Blk&\cr&&&\cr} 
+ (q^4+q^5) \Young{\Blk&\Blk&\Blk&\cr&&&&\cr} + q^6 \Young{&&&&&\cr}\cr
&= q^6 \Young{&&&&&\cr}}\mark$$

We have therefore computed that $H_{(2^3)}[X;q]=s_{(2^3)} + (q +q^2)s_{(321)}
+ q^3 s_{(33)} +
q^3 s_{(4,1,1)} + (q^2+q^3+q^4) s_{(4,2)}  +  (q^4+q^5 )s_{(51)} + q^6 s_{(6)}$.

\section{Hats and Ribbons}

Let $\Lambda$ represent the space of symmetric functions on an arbitrary
number of variables considered as the polynomials over $\Q$ in
the power symmetric functions $\{ p_1, p_2, p_3,\ldots \}$.
If $\X = \{x_1,x_2,x_3, \ldots\}$ is a set of variables then
denote the symmetric functions in these variables by $\Lambda^X$.  These
two spaces are isomorphic and here we will often identify the two.

We will make use of a mix of plethystic notation and
notation made standard by Macdonald in [3].
Plethystic notation is a device for expressing the substitution of
the monomials of one expression in a symmetric function.
Say that $E$ is formal series in a
set of variables $x_1, x_2, ...$ with possible special parameters
$q$ and $t$ which should be thought of as unknown elements of $\Q$.
For $k \geq 1$, define $p_k[E]$ to be $E$ with $x_i$ replaced by  
$x_i^k$ and $q$ and $t$ replaced by $q^k$ and $t^k$ respectively.
For a symmetric function $P$, $P[E]$ will represent the
the formal series found by expanding $P$ in terms
of the power symmetric functions and then substituting $p_k[E]$
for $p_k$.
More precisely, if the power sum expansion of the symmetric
function $P$ is given by
$$P = \sum_\la c_\la p_\la \mark$$  then $P[E]$ is given by the formula
$$P[E] = \sum_\la c_\la p_{\la_1}[E] p_{\la_2}[E] \cdots p_{\la_{\ell(\la)}}[E]. \mark$$

To evaluate a symmetric function in a set of variables $\left\{
x_1, x_2, x_3, \ldots \right\}$, set $X = x_1+x_2+x_3+ \cdots$ and then
we have that $p_k[X] = x_1^k + x_2^k + x_3^k + \cdots$ and
$P[X]$ represents the symmetric function $P$ evaluated at the $x_i$'s.
For this exposition we will use capital letters $X$ and $Y$ to represent 
sums of infinite sets of variables, $x_i$'s and $y_i$'s respectively.

There is a well known scalar product on $\Lambda$ defined by setting
$$\left< p_\lambda, p_\mu \right> = \delta_{\la\mu} z_\lambda 
= \delta_{\la\mu} \prod_{i \geq 1} i^{n_i(\la)}
n_i(\lambda)! \mark$$
where $\delta_{xy}$ is the Kronecker delta and 
we have used the notation $n_i(\la)$ to represent the number of parts
of size $i$ in $\la$.

For any symmetric function $f$, we denote the operation of skewing by $f$
by $f^\perp$ to represent the operation dual to multiplication by $f$ with
respect to this scalar product. More precisely, $\left< f^\perp g, h \right> =
\left< g, f h \right> $ and
$f^\perp g = \sum_\la s_\la \left< s_\la f, g \right>$.

For the Schur functions and the Hall-Littlewood symmetric functions, the
vertex operators that add a row to the indexing symmetric function are well
known.  Define the two symmetric function operators and their generating
functions $S(z)$ and $H(z)$ by
$$S_m P[X]  = P\left[ X- {1 \over z} \right] \Omega[zX] \coeff_{z^m} 
= S(z) P[X] \coeff_{z^m}\mark$$
$$H_m^q P[X]  = P\left[ X - {1-q \over z} \right] \Omega[zX] \coeff_{z^m} 
= H(z) P[X] \coeff_{z^m}
\mark$$
where $\Omega[X] = \sum_{n \geq 0} h_n[X]$. For a partition $\la$,
$S_{\la_1} S_{\la_2} \cdots S_{\la_\ell} 1 = s_\la[X]$ and
$H_{\la_1}^q H_{\la_2}^q \cdots H_{\la_\ell}^q 1 = H_\la[X;q] :=
\sum_{\mu \vdash |\la|} K_{\la\mu}(q) s_\mu[X]$ where $K_{\la\mu}(q)$
are the Kostka-Foulkes coefficients.  The operator $S_m$ is
due to Bernstein [5], and the operator $H_m^q$ is due to Jing [2]
although the notation and presentation here follows more closely [1].

Let $P[X]$ be an arbitrary symmetric function in the set of $X$ variables
and let $V$ be an operator (an element of $Hom(\Lambda^X, \Lambda^X)$).
Define $\ht{V}$ by
its action on $P[X]$ by the following formula

$$\ht{V} P[X] = (V^Y P[X-Y] ) \coeff_{Y=X} 
\mark$$
The $Y$ in the $V^Y$ is there to emphasize that $V$ is acting in the
`dummy' $Y$ set of variables only.  
After that operation is complete we set the $Y$ variables
equal to the $X$ variables.

The hat operation also appears in a study of operators that add rows and
columns to the standard bases of the symmetric functions [4].  The notation
used here makes it a useful tool for deriving identities.  The following 
proposition shows an important property of the hat operation, that it is 
an involution.

\Prop  Let $V$ be an element of $Hom(\Lambda,\Lambda)$. 
$\ht{\ht{V}} = V$.

\vskip .1in
\noindent
\Proof
Let $P[X]$ be an arbitrary symmetric function
$$\eqalign{\ht{\ht{V}} P[X] &= {\ht{V}}^Y P[X-Y] \coeff_{Y=X}\cr
&= V^Z P[X-(Y-Z)] \coeff_{Z=Y} \coeff_{Y=X}\cr
&= V^Z P[X-Y+Z] \coeff_{Y=X} \coeff_{Z=X} \cr
&= V^Z P[Z] \coeff_{Z=X} = V P[X]}\mark$$
\endofproof

If $R \models k$ and $R^+ \models k+1$ such that
$D(R) = D(R^+)$ (that is, $R^+$ is $R$ with
a cell to the right) then it follows directly from the
definition
$$S^{R} \Sw_1 = S^{R^+}.
\mark$$

A direct computation shows the following astonishing
relation with the hat
involution which provides a method for adding a 
cell bellow all others in the ribbon.

\Thm   If $R \models k$
and $R_+ \models k+1$ such that $D(R)\cup \{k\} = D(R_+)$, then
$$\ht{\ht{S^{R}} S_1} = S^{R_+}\mark$$

Before proceeding with the proof of this theorem, we remark that 
equations $(26)$ and $(27)$ imply that
the following recursive definition
is equivalent to the definition of $H_{1^k}^q$ given in equation $(8)$. 

Let
$H_{1^0}^q = \Sw_0$, $H_{1^1}^q = \Sw_1$, and set
$$H_{1^k}^q = q^{k-1} H_{1^{k-1}}^q \Sw_1 + \ht{\ht{H_{1^{k-1}}^q}S_1}\mark$$

We will also need the following two lemmas to prove Theorem $4$.

\Lem  For any operator $V$, 
$\ht{\ht{V}S_m} = \sum_{j\geq 0}(-1)^{m-j} h_j V \Sw_{m-j}$

\Proof  Let $V$ be any operator.
$$\eqalign{\ht{\ht{V} S_m }&= V^WP[X - ( Y-W) + 1/z] \Omega[z(Y-W)]
\coeff_{W=Y} \coeff_{z^m} \coeff_{Y=X}\cr
&= \Omega[zX] V^W P[W + 1/z] \Omega[-zW] \coeff_{W=X} \coeff_{z^m}}\mark$$
Set $\Sw(z) P[X] = P\left[X + {1 \over z }\right] \Omega[-zX]$, then 
$\Sw(z) P[X] \coeff_{z^\ell} = (-1)^\ell \Sw_\ell P[X]$.
Take the coefficient of $z^m$ to yield the identity.
\endofproof

\Lem
 $\ht{\ht{s_{\la}^\perp}S_{-m}} = 
(s_{(m,\la)})^\perp$.  

\Proof The operator $\ht{s_\la^\perp}$ has the action
$$\ht{s_\la^\perp} s_\mu[X] = s_\la^{\perp Y} s_\mu[X-Y] \coeff_{Y=X}
= \sum_{\gamma} (-1)^{|\gamma|} s_{\mu\slash \gamma}[X] s_{\gamma'\slash\la}[X]
=\left\{ {(-1)^{|\la|} \ if \ \la = \mu' \atop
0 \ \ \ otherwise }\right. \mark$$

Take the coefficient of $s_\mu[Y]$ in the following equation.
$$\eqalign{\ht{\ht{s_{\la}^\perp}S_{-m}}^X \Omega[XY] &= 
\ht{s_\la^\perp}^U S_{-m}^U \Omega[XY] \Omega[-YU] \coeff_{U=X} \cr
&= \Omega[XY] \ht{s_\la^\perp}^U S_{-m}^U \Omega\left[-Y\left(U-{1\over z}
\right)\right]
\Omega[zU] \coeff_{U=X}\cr
&=\Omega[XY] \ht{s_\la^\perp}^X \Omega[-(Y-z)X]\Omega\left[{Y\over z}\right]
\coeff_{z^{-m}}\cr
&=\Omega[XY] s_\la[Y-z]\Omega\left[{Y\over z}\right]
\coeff_{z^{-m}}\cr
&= \Omega[XY] s_{(m,\la)}[Y]} \mark$$
  On the left
one has $\ht{\ht{s_{\la}^\perp}S_{-m}} s_\mu[X]$
and on the right $s_{\mu \slash (m,\la)}[X]$.\endofproof

\Proof (of Theorem $4$)
Say that $S^R = s_{\la^{rc}}^\perp \Sw_{\la'}$ for some partition
$\la$ and for brevity we have used the notation $\Sw_{\la'} := \Sw_{\la_1'}
\Sw_{\la_2'} \cdots \Sw_{\la_m'}$ with $\la_1 = m$.

$$\eqalign{ \ht{\ht{S^{R}}S_1} &= 
\sum_{i\geq 0}(-1)^{1-i}  h_i s_{\la^{rc}}^\perp 
\Sw_{\la'} \Sw_{1-i} \cr
&= \sum_{i\geq 0}(-1)^{1+i+m}  h_i s_{\la^{rc}}^\perp 
\Sw_{1-i-m} \Sw_{(m,\la)'} \cr
&= \ht{\ht{s_{\la^{rc}}^\perp}S_{1-m}}\Sw_{(m,\la)'}}\mark$$
By the previous lemma this operator
reduces to $s_{(m-1,\la^{rc})}^\perp \Sw_{(m,\la)'}$, which is exactly
$S^{R_+}$
\endofproof

$(26)$ and $(27)$ provide a method for building all of the ribbon operators
recursively. It follows in our next theorem 
that the hat involution sends ribbon operators
to ribbon operators and permutes them in a very
natural and non-trivial manner.

Say that two operators are $k$-level equal and write $U \keq V$
if $V(s_\la) = U(s_\la)$ for
all $\ell(\la) \leq k$.

\Thm  Let $R$ be a ribbon of size $k$, then
$$\omega \ht{S^R} \omega \keq S^{\ht{R}}. \mark$$

We make the following general remark about operators and their hats before
proceeding with the proof.

\Lem
$U \keq V$ implies that $\omega \ht{U} \omega \keq \omega \ht{V} \omega$

\Proof 
If $\ell(\gamma) \leq k$ and $U(s_\la) = V(s_\la)$ for all $\la$ such that
$\ell(\la) \leq k$ then
$$\omega \ht{U} \omega (s_\gamma) 
= \sum_{\la \subseteq \gamma} (-1)^{|\la|} 
\omega ( U (s_\la) s_{\gamma' \slash \la'})
= \sum_{\la \subseteq \gamma} (-1)^{|\la|} 
\omega (V (s_\la) s_{\gamma' \slash \la'})
= \omega \ht{V} \omega (s_\gamma)\mark$$
\endofproof

\Proof (of Theorem $7$)
  The proof proceeds by induction on $k$.  The result is true for
$R \models 1$ since $\omega \ht{\Sw_1}\omega (s_{(n)})=s_{({n+1})}
= \Sw_1 (s_{(n)})$.

Let $R^+$ be a ribbon of size $k+1$.  Now either $k \in D(R^+)$ and $D(R^+) = 
D(R) \cup \{k\}$ for some ribbon $R$ of size $k$, or 
$k \notin D(R^+)$ and $D(R^+) = D(R)$ for some ribbon $R \models k$.

In the first case we have that $\omega \ht{S^{R^+}} \omega = \omega \ht{S^R}
S_1 \omega = \omega \ht{S^R}\omega \Sw_1$.  Note that if $\ell(\la) = k+1$ and
$\Sw_1(s_\la) = \pm s_\mu$ then 
$\ell(\mu) = \ell(\la) -1$. Since $\omega \ht{S^R}\omega$
is $k$-level equal to $S^{\ht{R}}$ then
we have  $\omega \ht{S^R}\omega \Sw_1 \simeq_{k+1} S^{\ht{R}} \Sw_1 
= S^{\ht{R^+}}$.

In the second case, the same reasoning implies $S^{R^+} = S^R \Sw_1 \simeq_{k+1}
\omega \ht{S^{\ht{R}}} \omega \Sw_1 = \omega\ht{S^{\ht{R}}}S_1 \omega$.
By Theorem $4$ and Lemma $8$, this implies $\omega \ht{S^{R^+}} \omega 
\simeq_{k+1} \ht{\ht{S^{\ht{R}}}S_1} 
= S^{\ht{R^+}}$. \endofproof

\section{Proof of Theorem $1$}

\Lem For $k \geq 2$
$$q^{k-1} H_{1^{k}}^q \Sw_2 = \ht{\ht{H_{1^{k}}^q} S_{2}}\mark$$

\Proof
The left hand side of this equation may be expanded using 
$\Sw_1 \Sw_2 = 0$, equation $(28)$ and Lemma $5$, we have  

$$\eqalign{q^{k-1} H_{1^{k}}^q \Sw_2 &= q^{k-1} \left( q^{k-1} H_{1^{k-1}}^q 
\Sw_1 + \ht{\ht{H_{1^{k-1}}^q}S_1}\right) \Sw_2 \cr
&= q^{k-1} \ht{\ht{H_{1^{k-1}}^q}S_1} \Sw_2  }\mark $$

The right hand side of the Lemma may also be expanded using the same
relations.  It follows that
$$\eqalign{\ht{\ht{H_{1^{k}}^q} S_{2}} &=
\sum_{j\geq 0} (-1)^{j+2} h_j H_{1^{k}}^q \Sw_{2-j}\cr
 &= \sum_{j\geq 0} (-1)^{j+2} h_j \left( q^{k-1} H_{1^{k-1}}^q \Sw_1 
+ \ht{\ht{H_{1^{k-1}}^q}S_1}\right)
 \Sw_{2-j}\cr
&= q^{k-1} \sum_{j\geq 0} (-1)^{j+2} h_j H_{1^{k-1}}^q \Sw_1 \Sw_{2-j} 
- \ht{\ht{H_{1^{k-1}}^q}S_1 S_2}\cr
&= q^{k-1} \ht{\ht{H_{1^{k-1}}^q}S_1} \Sw_2 }\mark$$ 
hence $(36)$ and $(37)$ are equal.\endofproof

\Proof {\sl (of Theorem $1$)}  Verify the following two
identities by a direct computation.
$$H(u) \Sw(z) P[X] = (q - u/z) \Sw(z) H(u) P[X]\mark$$
$$H(u) \Omega[zX] P[X] = {1-z/u \over 1-q z/u }\Omega[zX] H(u) P[X].\mark$$

Assume by induction for $\ell<k$, $H_{m+1}^q H_{1^\ell}^q = 
H_{1^{\ell+1}}^q H_{m}^q$.  Then in particular $H(u) H_{1^{k-1}}^q
= u H_{1^k}^q H(u)$.  By equation $(29)$, we have
$$\eqalign{ H_{m+1}^q \ht{\ht{H_{1^{k-1}}^q} S_1} &= 
H(u) \Omega[z X] H_{1^{k-1}}^q \Sw(z)  \coeff_{z^1} \coeff_{u^{m+1}}\cr
&={1-z/u \over 1-q z/u }\Omega[zX] u H_{1^{k}}^q (q - u/z) \Sw(z) H(u) 
\coeff_{z^1} \coeff_{u^{m+1}}\cr
&=({u-u^2/z})\Omega[zX] H_{1^{k}}^q \Sw(z) H(u) 
\coeff_{z^1} \coeff_{u^{m+1}} \cr 
&= \ht{\ht{H_{1^{k}}^q} S_1} H_m^q
- \ht{\ht{H_{1^{k}}^q} S_2} H_{m-1}^q} \mark$$
Using the same relations, we compute
$$\eqalign{ H_{m+1}^q H_{1^{k-1}}^q \Sw_1 &= 
H(u) H_{1^{k-1}}^q \Sw(z) \coeff_{z^1} \coeff_{u^{m+1}}\cr
&= u H_{1^{k}}^q (q - u/z) \Sw(z) H(u) \coeff_{z^1} \coeff_{u^{m+1}}\cr
&=({q u-u^2/z}) H_{1^{k}}^q \Sw(z) H(u)
\coeff_{z^1} \coeff_{u^{m+1}}\cr
&=q H_{1^{k}}^q \Sw_1 H_m^q+
H_{1^{k}}^q \Sw_2 H_{m-1}^q}\mark$$

Now using the recursive definition in equation $(28)$ for 
$H_{1^k}^q$ we have

$$H_{m+1}^q H_{1^k}^q =
\ht{\ht{H_{1^{k}}^q} S_1} H_m^q 
+ q^k H_{1^{k}}^q \Sw_1 H_m^q
- \ht{\ht{H_{1^{k}}^q} S_2} H_{m-1}^q
+ q^{k-1} H_{1^{k}}^q \Sw_2 H_{m-1}^q\mark$$

By Lemma $9$ and equation $(28)$, the right hand side reduces to 
$H_{1^{k+1}}^q H_m^q$.
\endofproof

\vskip .2in
\noindent
Acknowledgment:  Thank you to Adriano Garsia for his support
and valuable discussions in the development of this problem and
the simplification of the proofs.

\section{References}

\vskip .2in

\article $1$|A.~M.~Garsia|Orthogonality of Milne's polynomials and
raising operators|Discrete Math.|99|(1992)|247--264|

\article $2$|N.~Jing|Vertex operators and Hall-Littlewood symmetric
functions|Adv. Math.|87|(1991)|226--248|
 
\livre $3$|I.~G.~Macdonald|Symmetric Functions and Hall
Polynomials, Oxford Mathematical Monographs|second edition, 
Oxford Univ. Press, 1995|

\unarticle $4$|M.~Zabrocki|Vertex operators for standard bases of the
symmetric functions|J. of Alg. Comb.|

\livre $5$|A.~V.~Zelevinsky|Representations of finite classical groups:
a Hopf algebra approach|Springer Lecture Notes, 869, 1981|

\end